\documentclass[10pt,portrait,reqno,12pt]{amsart}%
\usepackage{amssymb,amsthm,amsmath, amsfonts, float}
\usepackage{graphicx}
\usepackage{multirow}
\usepackage{subfig}
\usepackage{amsmath}
\usepackage{amsfonts}
\usepackage{lscape}
\usepackage{amssymb}%
\providecommand{\U}[1]{\protect\rule{.1in}{.1in}}
\theoremstyle{plain}

\newtheorem{example}{Example}

\numberwithin{equation}{section}
\begin{document}
\title[Loxodromes on Twisted Surfaces in $E^{3}$]{Loxodromes on Twisted Surfaces in Euclidean 3-Space}
\subjclass[2010]{53A04, 53A05}
\keywords{Loxodromes, Meridian, Twisted surface.}
\author[M. Alt\i n]{\bfseries Mustafa Alt\i n$^{1\ast}$}
\address{$^{1}$Technical Sciences Vocational School, Bing\"{o}l University, Bing\"{o}l, Turkey\\
 \\
$^{\ast}$Corresponding author: maltin@bingol.edu.tr}

\begin{abstract}
In the present paper, loxodromes, which cut all meridians and parallels of
twisted surfaces (that can be considered as a generalization of rotational
surfaces) at a constant angle, have been studied in Euclidean 3-space and some
examples have been constructed to visualize and support our theory.

\end{abstract}
\maketitle


\section{General Information and Basic Concepts}

Loxodromes (also known as rhumb lines) correspond to the curves which
intersect all of the meridians at a constant angle on the Earth (see Figure
1). An aircraft flying and a ship sailing on a fixed magnetic compass course
move along a curve. Here the course is a rhumb and the curve is a loxodrome.
Generally, a loxodrome is not a great circle, thus it does not measure the
shortest distance between two points on the Earth. However loxodromes are
important in navigation and they should be known by aircraft pilots and
sailors \cite{Alex}.

If the shape of the Earth is approximated by a sphere, then the loxodrome is a
logarithmic spiral that cuts all meridians at the same angle and
asymptotically approaches the Earth's poles but never meets them. Since
maritime surface navigation defines the course as the angle between the
current meridian and the longitudinal direction of the ship, it may be
concluded that the loxodrome is the curve of the constant course, which means
that whenever navigating on an unchanging course we are navigating according
to a loxodrome \cite{Sergio}.

\begin{figure}[H]
\centering
\includegraphics[
    height=2.3in, width=2.8in
    ]{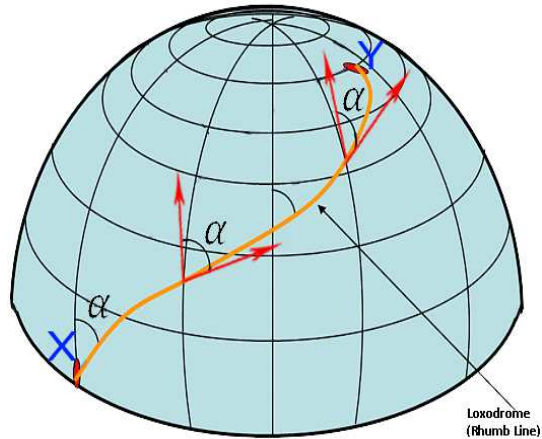}\caption{Loxodrome on Earth}%
\label{fig:y}%
\end{figure}

In this context, there are lots of studies about loxodromes in Euclidean and
Minkowskian spaces. For instance, the differential equations of loxodromes on
a sphere, spheroid, rotational surface, helicoidal surface and canal surface
in Euclidean 3-space have been given in \cite{Sergio}, \cite{Miljen},
\cite{Kos}, \cite{Baba1} and \cite{Baba2}, respectively. Also, in \cite{Baba3}
and \cite{Baba4}, spacelike and timelike loxodromes on rotational surface and
in \cite{Baba5}, differential equations of the spacelike loxodromes on the
helicoidal surfaces in Lorentz-Minkowski 3-space have been given.

Now, let we recall some basic notions about curves and twisted surfaces in
Euclidean 3-space $E^{3}$.

For two vectors $\vec{u}=(u_{1},u_{2},u_{3})$ and $\vec{v}=(v_{1},v_{2}%
,v_{3})$ in $E^{3}$, the inner product of these vectors and the norm of the
vector $\vec{u}$ are defined by%
\begin{equation}
\left\langle \vec{u},\vec{v}\right\rangle =u_{1}v_{1}+u_{2}v_{2}+u_{3}v_{3}
\label{1}%
\end{equation}
and%
\begin{equation}
\left\Vert \vec{u}\right\Vert =\sqrt{\left\langle \vec{u},\vec{u}\right\rangle
}, \label{2}%
\end{equation}
respectively. We say that $\vec{u}$ is a unit vector, if it satisfies
$\left\Vert \vec{u}\right\Vert =1$.

The arc-length of a regular curve $\alpha:I\subset%
\mathbb{R}
\longrightarrow E^{3},$ $s\longrightarrow$ $\alpha(s),$ between $s_{0}$ and
$s$ is%
\begin{equation}
t(s)=\int\nolimits_{s_{0}}^{s}\left\Vert \alpha^{\prime}(s)\right\Vert ds.
\label{3}%
\end{equation}
Then the parameter $t\in J\subset%
\mathbb{R}
$ is determined as $\left\Vert \alpha^{\prime}(t)\right\Vert =1.$

Also, the angle $\phi$ $(0<\phi<\pi)$ between the vectors $u$ and $v$ is%
\begin{equation}
\cos\phi=\frac{\left\langle \vec{u},\vec{v}\right\rangle }{\left\Vert \vec
{u}\right\Vert \left\Vert \vec{v}\right\Vert }. \label{4}%
\end{equation}

Here, we recall the definition and parametrization of twisted surfaces in
$E^{3}.$ (For detail, see \cite{Goemans}.)

A twisted surface in $E^{3}$ is obtained by rotating a planar curve $\alpha$
in its supporting plane while this plane itself is rotated about some
containing straight line. Without loss of generality, the coordinate system
can be chosen in such a way that the $xz$-plane corresponds with the plane
supporting the planar curve with the $z$-axis as its containing rotation axis
and that the straight line through the point $(a,0,0)$ parallel with the
$y$-axis acts as rotation axis for the planar curve.

Firstly, let we apply the rotation about the straight line through the point
$(a,0,0)$ parallel with the $y$-axis to the profile curve $\alpha
(y)=(f(y),0,g(y))$ ($f$ and $g$ are real-valued functions) and next apply the
rotation about the $z$-axis to the obtained surface. Then, up to a
transformation, we get the parametrization of the twisted surface in $E^{3}$
as%
\begin{equation}
T(x,y)=\left(
\begin{array}
[c]{c}%
(a+f(y)\cos(bx)-g(y)\sin(bx))\cos x,(a+f(y)\cos(bx)-g(y)\sin(bx))\sin x,\\
f(y)\sin(bx)+g(y)\cos(bx)
\end{array}
\right)  .\label{5}%
\end{equation}
Here, the presence of the factor $b\in%
\mathbb{R}
$ allows for differences in the rotation speed of both rotations and it is
obvious from the construction that, if we take $b=0,$ then the twisted surface
reduces to a surface of revolution. Thus, the twisted surfaces can be
considered as generalizations of surfaces of revolution.

After giving the definition of the twisted surfaces, twisted surfaces with
null rotation axis in Minkowski 3-space have been studied in \cite{Goemans2}
and twisted surfaces with vanishing curvature in Galilean 3-space have been
classified in \cite{Dede}. Also, in \cite{Kazan}, the twisted surfaces in
pseudo-Galilean space have been studied.

\section{Loxodromes on Twisted Surfaces in $E^{3}$}

In this section, we obtain the equations of loxodromes on the twisted surfaces
in $E^{3}.$

Let $T$ be the twisted surface which parametrized as (\ref{5}). Then the
coefficients of first fundamental form of the twisted surface $T$ are obtained
by%
\begin{equation}
\left.
\begin{array}
[c]{l}%
g_{11}=\frac{1}{2}\left\{
\begin{array}
[c]{l}%
2a^{2}+(1+2b^{2}+\cos(2bx))f^{2}+(1+2b^{2}-\cos(2bx))g^{2}\\
-4ag\sin(bx)+4f\cos(bx)(a-g\sin(bx))
\end{array}
\right\}  ,\\
g_{12}=g_{21}=b\left\{  fg^{\prime}-f^{\prime}g)\right\}  ,\\
g_{22}=f^{\prime2}+g^{\prime2},
\end{array}
\right\}  \label{6y}%
\end{equation}
where $g_{11}=\left\langle T_{x},T_{x}\right\rangle ,$ $g_{12}=g_{21}%
=\left\langle T_{x},T_{y}\right\rangle ,$ $g_{22}=\left\langle T_{y}%
,T_{y}\right\rangle ,$ $f=f(y),$ $g=g(y),$ $f^{\prime}=\frac{df(y)}{dy}$ and
$g^{\prime}=\frac{dg(y)}{dy}$.

Also, we know that, the first fundamental form in the base $\{M_{x},M_{y}\}$
for a surface $M(x,y)$ is given by%
\begin{equation}
ds^{2}=g_{11}dx^{2}+2g_{12}dxdy+g_{22}dy^{2},\label{7}%
\end{equation}
where $g_{ij}$ are the coefficients of the first fundamental form of $M$. So,
from (\ref{6y}) and (\ref{7}), we can write the first fundamental form of the
twisted surface (\ref{5}) as%
\begin{align}
ds^{2} &  =\left\{
\begin{array}
[c]{l}%
a^{2}-2ag\sin(bx)+2f\cos(bx)(a-g\sin(bx))\\
+\frac{1}{2}\left(  (1+2b^{2}+\cos(2bx))f^{2}+(1+2b^{2}-\cos(2bx))g^{2}%
\right)
\end{array}
\right\}  dx^{2}\nonumber\\
&  \text{ \ \ }+2b\left\{  fg^{\prime}-f^{\prime}g)\right\}  dxdy+(f^{\prime
2}+g^{\prime2})dy^{2}\label{8}%
\end{align}
and from (\ref{8}), the arc-length of any curve on the twisted surface between
$x_{1}$ and $x_{2}$ is given by%
\begin{equation}
s=\left\vert
{\displaystyle\int\nolimits_{x_{1}}^{x_{2}}}
\sqrt{%
\begin{array}
[c]{l}%
\left\{
\begin{array}
[c]{l}%
a^{2}-2ag\sin(bx)+2f\cos(bx)(a-g\sin(bx))\\
+\frac{1}{2}\left(  (1+2b^{2}+\cos(2bx))f^{2}+(1+2b^{2}-\cos(2bx))g^{2}%
\right)
\end{array}
\right\}  \\
+2b\left\{  fg^{\prime}-f^{\prime}g)\right\}  \frac{dy}{dx}+(f^{\prime
2}+g^{\prime2})\left(  \frac{dy}{dx}\right)  ^{2}%
\end{array}
}\text{ \ }dx\right\vert .\label{9}%
\end{equation}

Furthermore, a curve $\gamma$ is called a \textit{loxodrome} on the twisted
surface $T$ in $E^{3}$ if it cuts all meridians ($y$ constant) (or parallels
($x$ constant)) of $T$ at a constant angle.

Now, let us suppose that $\gamma(t)=T(x(t),y(t))$; i.e. $\gamma(t)$ is a curve
on the twisted surface $T$. With respect to the local base $\{T_{x},T_{y}\}$,
the vector $\gamma^{\prime}(t)$ has the coordinates $(x^{\prime},y^{\prime})$
and the vector $T_{x}$ has the coordinates $(1,0)$. At the point $p=T(x,y)$;
where the loxodrome cuts the meridians at a constant angle; we get%
\begin{equation}
\cos\phi=\frac{\left\langle \gamma^{\prime}(t),T_{x}\right\rangle }{\left\Vert
\gamma^{\prime}(t)\right\Vert \left\Vert T_{x}\right\Vert }=\frac
{g_{11}dx+g_{12}dy}{\sqrt{g_{11}^{2}dx^{2}+2g_{11}g_{12}dxdy+g_{11}%
g_{22}dy^{2}}}.\label{10}%
\end{equation}
Therefore, from (\ref{10}) we get%
\begin{equation}
g_{11}^{2}\sin^{2}\phi\text{ }dx^{2}+2g_{11}g_{12}\sin^{2}\phi\text{
}dxdy+(g_{12}^{2}-g_{11}g_{22}\cos^{2}\phi)\text{ }dy^{2}=0\label{11}%
\end{equation}
and so,%
\begin{equation}
\frac{dy}{dx}=\frac{-2g_{11}g_{12}\sin^{2}\phi\mp g_{11}\sqrt{g_{11}%
g_{22}-g_{12}^{2}}\sin(2\phi)}{2(g_{12}^{2}-g_{11}g_{22}\cos^{2}\phi
)}.\label{12'}%
\end{equation}

Hence from (\ref{6y}) and (\ref{12'}), the loxodrome on the twisted surface
(\ref{5}) must satisfy the following equation
\begin{equation}
\frac{dy}{dx}=-\frac{\left(
\begin{array}
[c]{l}%
\left(
\begin{array}
[c]{l}%
2b\sin\phi(fg^{\prime}-gf^{\prime})\\
\mp2\cos\phi\sqrt{%
\begin{array}
[c]{l}%
-b^{2}(gf^{\prime}-fg^{\prime})^{2}\\
+\left\{
\begin{array}
[c]{l}%
a^{2}-2ag\sin(bx)+2f\cos(bx)(a-g\sin(bx))\\
+\frac{1}{2}\left(
\begin{array}
[c]{l}%
(1+2b^{2}+\cos(2bx))f^{2}\\
+(1+2b^{2}-\cos(2bx))g^{2}%
\end{array}
\right)
\end{array}
\right\}  (f^{\prime2}+g^{\prime2})
\end{array}
}%
\end{array}
\right)  \times\\
\\
\left(
\begin{array}
[c]{l}%
2a^{2}+(1+2b^{2}+\cos(2bx))f^{2}+(1+2b^{2}-\cos(2bx))g^{2}\\
-4ag\sin(bx)+4f\cos(bx)(a-g\sin(bx))
\end{array}
\right)  \sin\phi
\end{array}
\right)  }{4\left(
\begin{array}
[c]{l}%
b^{2}\left(  gf^{\prime}-fg^{\prime}\right)  ^{2}\\
-\left\{
\begin{array}
[c]{l}%
a^{2}-2ag\sin(bx)+2f\cos(bx)(a-g\sin(bx))\\
+\frac{1}{2}\left(  (1+2b^{2}+\cos(2bx))f^{2}+(1+2b^{2}-\cos(2bx))g^{2}%
\right)
\end{array}
\right\}  (f^{\prime2}+g^{\prime2})\cos^{2}\phi
\end{array}
\right)  }.\label{12}%
\end{equation}
Here, let us construct two examples, with the aid of Mathematica, to visualize
and support our theory.

\begin{example}
Taking the profile curve as $\alpha(y)=(y,0,0)$, the twisted surface (\ref{5})
becomes%
\begin{equation}
T(x,y)=\left(  (a+y\cos(bx))\cos x,(a+y\cos(bx))\sin x,y\sin(bx)\right)
.\label{13}%
\end{equation}
From (\ref{12}) (we take $\mp$ in this equation as $-$), we get the
differential equation of the loxodrome on the twisted surface (\ref{13}) for
$a=0$ and $b=\frac{1}{2}$ as%
\begin{equation}
\frac{dy}{dx}=\frac{y\sqrt{2\cos x+3}\tan\phi}{2}.\label{14}%
\end{equation}
So, we have%
\[
\frac{dy}{y}=\frac{\sqrt{2\cos x+3}\tan\phi}{2}dx
\]
and by integrating both sides of this equation, we get%
\begin{equation}
\ln y=%
{\displaystyle\int\nolimits_{x_{0}}^{x}}
\frac{\sqrt{2\cos x+3}\tan\phi}{2}dx.\label{15}%
\end{equation}
Putting $x_{0}=0$ in (\ref{15}), we reach that%
\begin{equation}
y=y(x)=e^{\sqrt{5}EllipticE\left[  \frac{x}{2},\frac{4}{5}\right]  \tan\phi
}.\label{16}%
\end{equation}
Now, if we take $\phi=\frac{\pi}{6}$ and $x\in(-2\pi,2\pi)$, we get
$y\in(0.0476989,20.9649).$ Thus, the loxodrome which lies on the twisted
surface (\ref{13}) is obtained as
\begin{equation}
\gamma(x)=e^{\sqrt{\frac{5}{3}}EllipticE\left[  \frac{x}{2},\frac{4}%
{5}\right]  }.\left(  \cos\left(  \frac{x}{2}\right)  \cos x,\cos\left(
\frac{x}{2}\right)  \sin x,\sin\left(  \frac{x}{2}\right)  \right)
.\label{17}%
\end{equation}
Also, the arc-length of our loxodrome (\ref{17}) is approximately equal to
$41.8343$. The twisted surface (\ref{13}), meridian for $y=15$ and the
loxodrome (\ref{17}) can be seen in Figure 2.
\end{example}

\begin{figure}[H]
\centering
\includegraphics[
    height=2.9in, width=3.5in
    ]{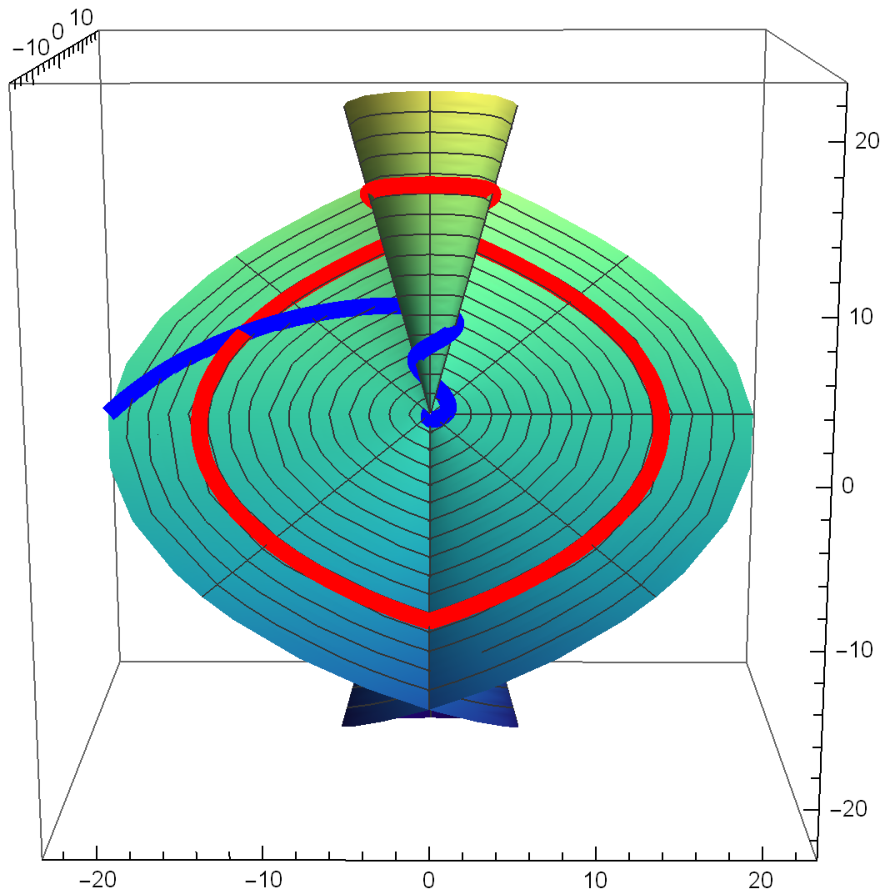}\caption{Twisted surface (\ref{13}), Meridian (red) and Loxodrome
(blue)}%
\label{fig:1}%
\end{figure}

\begin{example}
For the profile curve $\alpha(y)=(\cos y,0,\sin y)$, the twisted surface
(\ref{5}) is%
\begin{equation}
T(x,y)=\left(
\begin{array}
[c]{c}%
(a+\cos y\cos(bx)-\sin y\sin(bx))\cos x,(a+\cos y\cos(bx)-\sin y\sin(bx))\sin
x,\\
\cos y\sin(bx)+\sin y\cos(bx)
\end{array}
\right)  .\label{18}%
\end{equation}
If we take $a=1$ and $b=0$, then from (\ref{12}) (we take $\mp$ in this
equation as $-$), we have%
\begin{equation}
\frac{dy}{dx}=2\cos^{2}\left(  \frac{y}{2}\right)  \tan\phi.\label{19}%
\end{equation}
So, we have%
\[
\frac{dy}{2\cos^{2}\left(  \frac{y}{2}\right)  }=\tan\phi dx
\]
and by integrating both sides of this equation, we get%
\begin{equation}
\tan\left(  \frac{y}{2}\right)  =%
{\displaystyle\int\nolimits_{x_{0}}^{x}}
\tan\phi\text{ }dx.\label{20}%
\end{equation}
Taking $x_{0}=0$ in (\ref{20}), we reach that%
\begin{equation}
y=y(x)=2\left(  \arctan\left(  x\tan\phi\right)  +c\pi\right)  ,\text{ }c\in%
\mathbb{Z}
.\label{21}%
\end{equation}
Here, if we take $c=0,$ $\phi=\frac{\pi}{4}$ and $x\in(-\pi,\pi)$, we have
$y\in(-2.52525,2.52525).$ Thus, the loxodrome which lies on the twisted
surface (\ref{18}) is obtained as
\begin{equation}
\gamma(x)=\left(  ((1+\cos(2\arctan x))\cos x,(1+\cos(2\arctan x)\sin
x,\sin(2\arctan x)\right)  .\label{22}%
\end{equation}
Also, the arc-length of our loxodrome (\ref{22}) is approximately equal to
$7.1425$. One can see the twisted surface (\ref{18}), meridian for $y=1$ and
the loxodrome (\ref{22}) in Figure 3.
\end{example}

\begin{figure}[H]
\centering
\includegraphics[
    height=2.5in, width=4.3in
    ]{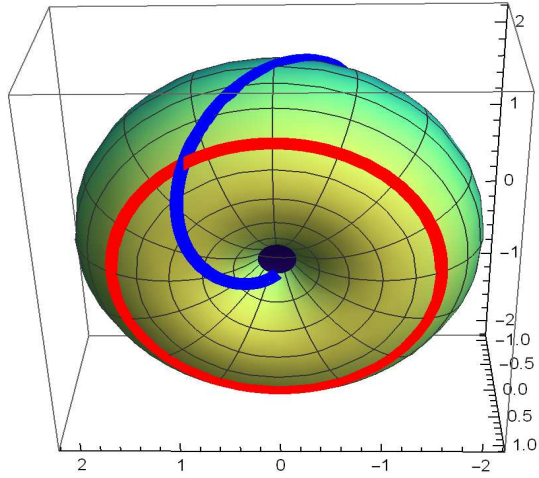}\caption{Twisted surface (\ref{18}), Meridian (red) and Loxodrome
(blue)}%
\label{fig:2}%
\end{figure}

Furthermore, from the definition of the angle $\theta$ between the loxodrome
and any parallel ($x$=constant), we have%
\begin{equation}
\cos\theta=\frac{\left\langle \gamma^{\prime}(t),T_{y}\right\rangle
}{\left\Vert \gamma^{\prime}(t)\right\Vert \left\Vert T_{y}\right\Vert }%
=\frac{g_{12}dx+g_{22}dy}{\sqrt{g_{11}g_{22}dx^{2}+2g_{12}g_{22}%
dxdy+g_{22}^{2}dy^{2}}}.\label{23}%
\end{equation}
From (\ref{23}), we get%
\begin{equation}
(g_{11}g_{22}\cos^{2}\theta-g_{12}^{2})\text{ }dx^{2}-2g_{12}g_{22}\sin
^{2}\theta\text{ }dxdy-g_{22}^{2}\sin^{2}\theta\text{ }dy^{2}=0\label{24}%
\end{equation}
and so,%
\begin{equation}
\frac{dx}{dy}=\frac{-2g_{12}g_{22}\sin^{2}\theta\mp g_{22}\sqrt{g_{11}%
g_{22}-g_{12}^{2}}\sin(2\theta)}{2(g_{12}^{2}-g_{11}g_{22}\cos^{2}\theta
)}.\label{24'}%
\end{equation}
Therefore, the loxodrome on the twisted surface must satisfy the following
equation
\begin{equation}
\frac{dx}{dy}=\frac{\left(
\begin{array}
[c]{l}%
2b\sin\theta(f^{\prime}g-fg^{\prime})\\
\mp2\cos\theta\sqrt{%
\begin{array}
[c]{l}%
-b^{2}(gf^{\prime}-fg^{\prime})^{2}\\
+\left\{
\begin{array}
[c]{l}%
a^{2}-2ag\sin(bx)\\
+2f\cos(bx)(a-g\sin(bx))\\
+\frac{1}{2}\left(
\begin{array}
[c]{l}%
(1+2b^{2}+\cos(2bx))f^{2}\\
+(1+2b^{2}-\cos(2bx))g^{2}%
\end{array}
\right)
\end{array}
\right\}  (f^{\prime2}+g^{\prime2})
\end{array}
}%
\end{array}
\right)  (f^{\prime2}+g^{\prime2})\sin\theta}{\left(
\begin{array}
[c]{l}%
2b^{2}\left(  gf^{\prime}-fg^{\prime}\right)  ^{2}\\
-\left(
\begin{array}
[c]{l}%
2a^{2}-4ag\sin(bx)+4f\cos(bx)(a-g\sin(bx))\\
+(1+2b^{2}+\cos(2bx))f^{2}+(1+2b^{2}-\cos(2bx))g^{2}%
\end{array}
\right)  (f^{\prime2}+g^{\prime2})\cos^{2}\theta
\end{array}
\right)  }\label{25}%
\end{equation}
Now, let us give an example for the loxodrome which cuts the parallels of the
twisted surface at a constant angle.

\begin{example}
Let us take the profile curve as $\alpha(y)=(\cos^{2}y,0,\sin^{2}y)$. Then,
the twisted surface (\ref{5}) is%
\begin{equation}
T(x,y)=\left(
\begin{array}
[c]{c}%
(a+\cos^{2}y\cos(bx)-\sin^{2}y\sin(bx))\cos x,(a+\cos^{2}y\cos(bx)-\sin
^{2}y\sin(bx))\sin x,\\
\cos^{2}y\sin(bx)+\sin^{2}y\cos(bx)
\end{array}
\right)  . \label{26}%
\end{equation}

Putting $a=-1$ and $b=0$, from (\ref{25}) (we take $\mp$ in this equation as
$-$), we get%
\begin{equation}
\frac{dx}{dy}=2\sqrt{2}\cot y\tan\theta.\label{27}%
\end{equation}
Thus, by integrating both sides of the equation%
\[
dx=2\sqrt{2}\cot y\tan\theta dy
\]
we have%
\begin{equation}
x=%
{\displaystyle\int\nolimits_{y_{0}}^{y}}
2\sqrt{2}\cot y\tan\theta dy.\label{28}%
\end{equation}
For $y_{0}=\frac{\pi}{2}$, we reach that%
\begin{equation}
x=x(y)=2\sqrt{2}\ln(\sin y)\tan\theta.\label{29}%
\end{equation}
Here, by taking $\theta=\frac{\pi}{3}$ and $y\in(\frac{\pi}{16},\frac{2\pi}%
{3})$, we have $x\in(-8.00637,-0.704674).$ Therefore, the loxodrome which lies
on the twisted surface (\ref{26}) is obtained as
\begin{equation}
\delta(y)=\left(  -\sin^{2}y\cos(2\sqrt{6}\ln(\sin y),-\sin^{2}y\sin(2\sqrt
{6}\ln(\sin y),\sin^{2}y\right)  .\label{30}%
\end{equation}
Also, the arc-length of the loxodrome (\ref{30}) is approximately equal to
$3.42788$. One can see the twisted surface (\ref{26}), meridian for $x=-1$ and
the loxodrome (\ref{30}) in Figure 4.
\end{example}

\begin{figure}[H]
\centering
\includegraphics[
    height=3.1in, width=3.5in
    ]{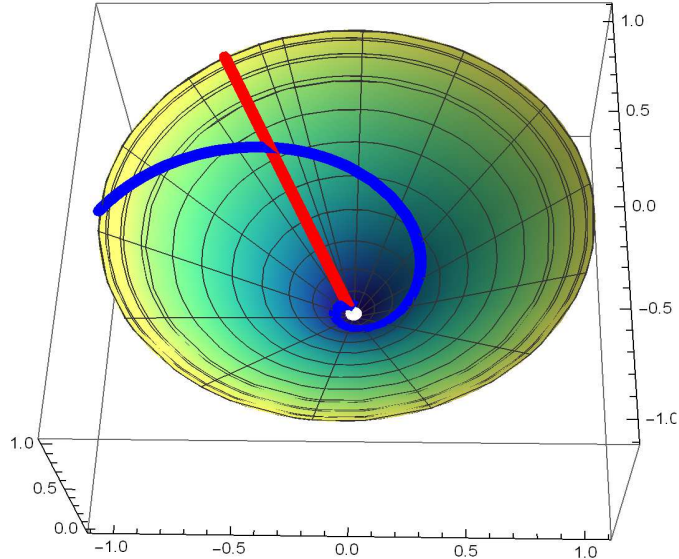}\caption{Twisted surface (\ref{26}), Meridian (red) and Loxodrome
(blue)}%
\label{fig:3}%
\end{figure}


\bigskip

\begin{thebibliography}{99}                                                                                               %


\bibitem {Alex}J. Alexander; \textit{Loxodromes: a rhumb way to go}, Math.
Mag.\textit{,} 77(5), (2004), 349--356.

\bibitem {Baba1}M. Babaarslan and Y. Yayl\i; \textit{Differential Equation of
the Loxodrome on a Helicoidal Surface, }The Journal of Navigation, 68, (2015), 962--970.

\bibitem {Baba2}M. Babaarslan, \textit{Loxodromes on Canal Surfaces in
Euclidean 3-Space}, Ann. Sofia Univ. Fac. Math and Inf., 103, (2016), 97--103.

\bibitem {Baba3}M. Babaarslan and Y. Yayl\i; \textit{Space-like loxodromes on
rotational surfaces in Minkowski 3-space}, J. Math. Anal. Appl., 409, (2014), 288--298.

\bibitem {Baba4}M. Babaarslan and M.I. Munteanu; \textit{Timelike loxodromes
on rotational surfaces in Minkowski 3--space}, Annals of the Alexandru Ioan
Cuza University-Mathematics, (2015), DOI: 10.2478/aicu-2013-0021.

\bibitem {Baba5}M. Babaarslan and M. Kayacik; \textit{Differential Equations
of the Space-Like Loxodromes on the Helicoidal Surfaces}, Differ Equ Dyn Syst,
28(2), (2020), 495--512.

\bibitem {Dede}M. Dede, C. Ekici, W. Goemans and Y \"{U}nl\"{u}t\"{u}rk;
\textit{Twisted surfaces with vanishing curvature in Galilean 3-space},
International Journal of Geometric Methods in Modern Physics, 15(1), (2018).

\bibitem {Goemans}W. Goemans and I. Van de Woestyne, \textit{Twisted surfaces
in Euclidean and Minkowski 3-space}, Pure and Applied Differential Geometry,
Joeri Van der Veken, Ignace Van de Woestyne, Leopold Verstraelen and Luc
Vrancken (Editors), Shaker Verlag (Aachen, Germany), (2013), 143-151.

\bibitem {Goemans2}W. Goemans and I. Van de Woestyne; \textit{Twisted Surfaces
with Null Rotation Axis in Minkowski 3-Space}, Results in Mathematics, 70,
(2016), 81--93.

\bibitem {Kazan}A. Kazan and H.B. Karada\u{g};\textit{ Twisted Surfaces in the
Pseudo-Galilean Space}, NTMSCI, 5(4), (2017), 72-79.

\bibitem {Sergio}S. Kos, D. Vranic and D. Zec;\textit{ Differential Equation
of a Loxodrome on a Sphere, }The Journal of Navigation, 52, (1999), 418--420.

\bibitem {Kos}S. Kos, R. Filjar and M. Hess; \textit{Differential equation of
the loxodrome on a rotational surface}, ION ITM Conference, (2009), Anaheim,
California, USA.

\bibitem {Miljen}M. Petrovic; \textit{Differential Equation of a Loxodrome on
the Spheroid}, "Na\v{s}e more", 54(3-4), (2007), 87-89.
\end{thebibliography}
\end{document}